\documentclass[12pt]{amsart}

\numberwithin{equation}{section}

\newcommand {\Z}{\mathbb{Z}}	
\newcommand {\C}{\mathbb{C}}  

\begingroup

\theoremstyle{plain}
\newtheorem{thm}{Theorem}[section]
\newtheorem{prop}[thm]{Proposition}

\newtheorem{cor}[thm]{Corollary}

\theoremstyle{definition}

\endgroup

\newcommand{\1}{{\bf 1}}

\newcommand{\myqed}{\ \hfill \rule{2mm}{2mm} \medskip }
\newcommand{\cd}{\cdots}

\newcommand{\noi}{\noindent}
\newcommand{\beq}{\begin{equation}}
\newcommand{\eeq}{\end{equation}}
\newcommand{\beqr}{\begin{eqnarray*}}
\newcommand{\eeqr}{\end{eqnarray*}}
\newcommand{\bc}{\begin{center}}
\newcommand{\ec}{\end{center}}

\newcommand{\ep}{\varepsilon}

\begin{document}

\title[Foguel-Hankel operators]{Schur multipliers and
operator-valued Foguel-Hankel operators}
\author{C. Badea}
\address{Math\'ematiques, UMR 8524 au CNRS, Universit\'e de Lille I,
  F--59655 Villeneuve d'Ascq, France}
\email{catalin.badea@agat.univ-lille1.fr}
\author{V.I. Paulsen}
\address{Department of Mathematics, University of Houston,
              Houston, TX 77204-3476}
\email{vern@math.uh.edu}
\date{}

\begin{abstract}
We show that some matrices are Schur multipliers and this is applied 
to
obtain classes of operator-valued Foguel-Hankel operators similar to
contractions. This provides partial answers to a problem of K. 
Davidson
and the second author concerning CAR-valued Foguel-Hankel operators.
\end{abstract}
\maketitle
\section{Introduction}
An example of a polynomially bounded operator on Hilbert space
not similar to a contraction
was found recently
by Pisier \cite{pisier}. An operator-theoretic proof
that certain CAR-valued
Foguel-Hankel operators are polynomially bounded operators but not 
similar to
contractions was given by Davidson and Paulsen
\cite{davidson/paulsen}. It is still an open question 
\cite{davidson/paulsen}
to characterize operators in this family which are similar to 
contractions.

The aim of this note is to prove some partial results concerning
this open problem. The present note is a sequel of
\cite{davidson/paulsen} where this problem is studied.
A certain
familiarity with \cite{davidson/paulsen} is supposed.
For the convenience of
the reader some
notation and known facts are recalled below.

\subsection{Background.} We denote by
$H$ a separable Hilbert space and by $B(H)$ the C*-algebra of all
bounded and linear operators on $H$. An operator $T\in B(H)$ is said
to be \emph{power bounded} if
\begin{equation}
	\sup_{n\in \Z_{+}} \|T^n\| < + \infty
	\label{eq:11}
\end{equation}
and \emph{polynomially bounded} if there exists a constant $K$ such
that
\begin{equation}
	\|p(T)\| \leq K\|p\|_{\infty}
	\label{eq:12}
\end{equation}
for each analytic polynomial $p$. Here 
$$\|p\|_{\infty} := \sup \{|p(z)| : |z| \leq 1\}.$$
We say that $T$ is  \emph{similar to a
contraction} if there is an invertible operator $L \in B(H)$ such that
$$
	\|L^{-1}TL\| \leq 1 .
$$

The following implications hold~:

\begin{eqnarray}
		T \mbox{ similar to a contraction } & \Rightarrow & T
\mbox{ polynomially
	bounded }
		\nonumber  \\
		 & \Rightarrow & T \mbox{ power bounded } .
	\label{eq:14}
\end{eqnarray}

The first implication follows from von Neumann's \cite{vNeumann}
inequality
$$
 	\|p(C)\| \leq \|p\|_{\infty},
$$
valid for each contraction $C \in B(H)$.
The second implication is clear from (\ref{eq:11}) and
 (\ref{eq:12}).

It was proved by Paulsen \cite{paulsen:jfa} that $T\in B(H)$ is 
similar to a
contraction if and only if $T$ is completely polynomially bounded,
that is, there exists a constant $K$ such that
$$
	\| \left[ p_{ij}(T) \right]_{1\leq i,j\leq n}\|
		\leq K \sup \{\|\left[ p_{ij}(z) \right]_{1\leq i,j\leq
n}\| : |z| \leq
		1\},
$$
for all positive integers $n$ and all $n\times n$ matrices $\left[
p_{ij}\right]_{1\leq i,j\leq n}$ with polynomial entries. Recall
that $\left[ p_{ij}(T)\right]_{1\leq i,j\leq n}$ is identified with
an operator acting on the direct sum of $n$ copies of the
corresponding Hilbert space in a natural way.

No implication in (\ref{eq:14}) can be reversed. The first power
bounded operator not similar to a contraction was constructed by
Foguel \cite{foguel}. His counterexample has the form
\begin{equation}
	R(X) = R(S^*,S ; X) =
	\left[\begin{array}{cc}
			S^* & X  \\
			0 & S
		\end{array}\right] ,
	\label{eq:16}
\end{equation}
where $S$ denotes the unilateral shift on $\ell_{2}$ and $X$ was a
suitable
diagonal projection onto a subspace of $\ell_{2}$. We will call
 \emph{Foguel operators} the
operators of type (\ref{eq:16}).

Lebow \cite{lebow} proved that the example constructed by Foguel is
not polynomially bounded. Other examples of power bounded, not
polynomially bounded operators are in \cite{davie, peller2,bozejko}.

A \emph{Foguel-Hankel operator} is a Foguel operator (\ref{eq:16})
with
$$X = [a_{i+j}]_{i,j\geq 0} = \Gamma_{f},$$
a Hankel operator with symbol $f$ (cf. \cite{davidson/paulsen}).
The study of Foguel-Hankel operators was initiated by Foias and 
Williams
\cite{foias/williams} and Peller \cite{peller}. It follows from the
work of Peller \cite{peller},
Bourgain \cite{bourgain} and Aleksandrov and
Peller \cite{aleksandrov/peller} that these Foguel-Hankel operators
are similar to contractions whenever
they are polynomially bounded. Both conditions are equivalent to
$f'\in \mbox{ BMOA }$, that is with the boundedness of
$$\Gamma_{f'} = [(i+j+1)a_{i+j}]_{i,j\geq 0}.$$

However, it is still unknown if a general Foguel operator is similar 
to a contraction whenever it is polynomially bounded. 
See \cite{ferguson}, where this is related to computing a certain Ext 
group.

The first example of a polynomially bounded operator not similar to a
contraction was found by Pisier \cite{pisier}. His
example is a CAR-valued Foguel-Hankel operator which we introduce 
below.

\subsection{CAR-valued Foguel-Hankel operators.} Let
$\Lambda$ be a function from $H$ into
$B(H)$ satisfying
the CAR - \textit{canonical anticommutation relations}~:
 for all $u,v\in H$,
\begin{gather*}\label{CAR}
 \Lambda(u)\Lambda(v)+\Lambda(v)\Lambda(u)=0  \\
\intertext{and}
 \Lambda(u)\Lambda(v)^*+\Lambda(v)^*\Lambda(u)= (u,v) I .
\end{gather*}
The range  of $\Lambda$ is isometric to Hilbert space.
Let $\{e_{n}\}_{n\geq 0}$ be an orthonormal basis for $H$, and let
$C_n=\Lambda(e_n)$ for $n\geq 0$. For an
arbitrary sequence $\alpha=(\alpha_0,\alpha_1,\dots)$ in $\ell^2$, let
$$
 Y_\alpha = \Bigl[ \alpha_{i+j}C_{i+j} \Bigr]
$$
be a CAR-valued Hankel operator. Let
$$
R(Y_{\alpha}) = R(S^{*(\infty)},S^{(\infty)} ; Y_\alpha) =
   \begin{bmatrix}S^{*(\infty)}&Y_\alpha\\0&S^{(\infty)}\end{bmatrix}
$$
be the
corresponding  \emph{CAR-valued Foguel-Hankel operator} \cite{pisier},
\cite{davidson/paulsen}.

The initial choice of
$\alpha$ made by Pisier was $\alpha_{2^k-1}=1$ for $k\ge0$ and
$\alpha_i=0$ otherwise. In this case $R(Y_\alpha)$ is polynomially
bounded but not completely polynomially bounded.

The following more
general result holds (cf.~\cite{pisier}, \cite{davidson/paulsen}).

For a fixed sequence $\alpha=(\alpha_0,\alpha_1,\dots) \in \ell_{2}$,
let
$$A = A(\alpha) := \sup_{k\ge0} (k+1)^2 \sum_{i\ge k}|\alpha_i|^2$$
and
$$B_2 = B_{2}(\alpha) := \sum_{k\ge0} (k+1)^2 |\alpha_k|^2 .$$
The operator $R(Y_{\alpha})$ is polynomially bounded if and only if
$A$
is finite. If $R(Y_{\alpha})$ is similar
to a
contraction, then $B_2$ is finite.

It is an open problem \cite{davidson/paulsen} to characterize in
terms of the sequence $\alpha$ when $R(Y_{\alpha})$ is similar
to a
contraction. In particular, it is not known
if $B_{2}(\alpha)$ finite implies the similarity of $R(Y_{\alpha})$
to a contraction. Note \cite[p. 163]{davidson/paulsen} that
$$B_{2}(\alpha) = \|\Gamma_{F'}\|^2 = 
\|(i+j+1)\alpha_{i+j}C_{i+j}\|^2,$$
where the operator-valued symbol $F$ given by
$$F(z) = \sum_{n\geq 0}\alpha_{n}C_{n}z^{-n-1}$$
is such that $\Gamma_{F} = Y_{\alpha}$.

We refer to \cite{pisier}, \cite{davidson/paulsen}, 
\cite{davidson:survey}
for more information and for the undefined terms.

\subsection{Organization of the paper.} The main results are
stated in the next section. Section three contains
some useful results about Schur multipliers. In the fourth section, a
sufficient condition for similarity is given.
These results are used
to prove the main results in the last section.

\section{Main results}
We use
notations as above. The first two results give sufficient conditions 
for
similarity to a contraction of an operator-valued Foguel
operator. Although these results are implicit in the work of 
\cite{foias/williams} and \cite{davidson/paulsen}, 
they do not seem to have been stated elsewhere.

\begin{thm}\label{thm:C}
	Let $X \in B(\ell_{2}(H))$ with matrix $X = [X_{ij}]_{i,j\geq 0}$,
	$X_{ij}\in
	B(H)$, with respect to a fixed orthonormal basis. For each $n \geq
	1$ set
	$$A^{(n)}(X)_{ij} = X_{ij} + X_{i-1,j+1} + \ldots
	X_{i-\min(i,n-1),j+\min(i,n-1)} $$
	and let $A^{(n)}(X)$ be the matrix $[A_{ij}^{(n)}(X)_{ij}]_{i,j\geq
0}$.
	If
\begin{equation}
	\sup_{n\geq 1}\|A^{(n)}(X)\| < +\infty,
	\label{eq:thmC}
\end{equation}
	then the operator-valued Foguel operator
	$$R(S^{*(\infty)}, S^{(\infty)} ; X) =
	\begin{bmatrix}S^{*(\infty)}&X\\0&S^{(\infty)}\end{bmatrix} \in
	B(\ell_{2}(H)\oplus \ell_{2}(H))$$
	is similar to a contraction.
\end{thm}

In the case of operator-valued Foguel-Hankel operators the following holds.

\begin{thm}\label{thm:charact}
    Let $\Gamma = [\Gamma_{i+j}]_{i,j\geq 0}$
be an operator-valued Hankel operator,
$\Gamma_{k}\in B(H)$. The operator-valued Foguel-Hankel operator
	$$R(S^{*(\infty)}, S^{(\infty)} ; \Gamma) =
	\begin{bmatrix}S^{*(\infty)}&\Gamma\\0&S^{(\infty)}\end{bmatrix} \in
	B(\ell_{2}(H)\oplus \ell_{2}(H))$$
	is similar to a contraction if any of the following operators 
	$$\Gamma D - D^*\Gamma = [(j-i)\Gamma_{i+j-1}]_{i,j\geq
0},$$
$$\Gamma D = [j \Gamma_{i+j-1}]_{i,j\geq 0},$$  or
$$D^*\Gamma = [i\Gamma_{i+j-1}]_{i,j\geq 0}$$
	is a bounded operator. Here $\Gamma_{-1} = 0$ and $D = 
	[(i+1)\delta_{i+1}^j]_{i,j\geq 0}$ is the
	differentiation operator.
\end{thm}

If $\Gamma = \Gamma_{f} = (a_{i+j})$ is a scalar Hankel operator with 
symbol $f$,
then boundedness of any one of the three operators above is 
equivalent to 
boundedness of the other two and this occurs if and only if $f'
\in \mbox{ BMOA}$. This is a consequence of the fact that both
conditions are equivalent to similarity to a contraction of the 
Foguel-Hankel
operator. This, as was remarked in \cite{aleksandrov/peller},
also follows from \cite{janson/peetre}.

For operator-valued Hankels, the situation is much more complicated 
as is shown in \cite{davidson/paulsen}.
A sufficient condition for similarity to a contraction for an
operator-valued Foguel-Hankel operator was given by Blower
\cite{blower} in terms of Carleson measures.

In the case of CAR-valued Foguel-Hankel operators, we
still do not know if $B_{2}(\alpha) < +\infty$ implies the similarity
to a contraction of $R(Y_{\alpha})$. Theorem \ref{thm:charact} and
the Schur multipliers results of the next section will imply the
following results.

\begin{thm}[$\log \log$ condition]\label{thm:A}
	Let $\ep > 0$. Suppose
	$$\sum_{k\geq 1} (k+1)^2\left[ \log
	(k+1)\right]^{2} \left[ \log (\log
	(k+1))\right]^{2+\ep}|\alpha_{k}|^2 < + \infty.$$
	Then the CAR-valued Foguel-Hankel operator $R(Y_{\alpha}) =
	R(S^{*(\infty)},S^{(\infty)} ;
	Y_{\alpha})$ is similar to a contraction.
\end{thm}

Since the logarithm goes to infinity less quickly than any power, we
obtain the following consequences.

\begin{cor}[$\log$ condition]\label{log}
Let $\ep > 0$. Suppose
	$$\sum_{k\geq 0} (k+1)^2\left[ \log
	(k+1)\right]^{2+\ep}|\alpha_{k}|^2 < + \infty.$$
	Then the CAR-valued Foguel-Hankel operator
	$R(Y_{\alpha})$ is similar to a contraction.
\end{cor}
and
\begin{cor}[power condition]\label{cor}
Let $\ep > 0$. Suppose
	$$B_{2+\ep}(\alpha)
	:= \sum_{k\geq 0} (k+1)^{2+\ep}|\alpha_{k}|^2 < +
	\infty.$$
	Then the CAR-valued Foguel-Hankel operator
	$R(Y_{\alpha})$ is similar to a contraction.
\end{cor}

A proof of Corollary \ref{cor}
can be given by combining results from
\cite{davidson/paulsen} and \cite{petrovic}.
A different proof in the case $\ep = 1$ can be
found in \cite{badea}.

\section{Schur multipliers}
Let $A = [a_{ij}]_{i,j\geq 1}$ and $B = [b_{ij}]_{i,j\geq 1}$
be two matrices of the same
size (finite or infinite). The Schur product of $A$ and $B$ is defined
to be the matrix of elementwise products
$$A\ast B = [a_{ij}b_{ij}]_{i,j\geq 1}.$$

For $M\in B(\ell_{2})$ we let $\mathcal{S}_{M}$ denote the Schur
multiplication map by $M$ on $B(\ell_{2})$, $\mathcal{S}_{M}(A) =
M\ast A$.
Then $M$ is said to be a \emph{Schur multiplier} if
$$\|\mathcal{S}_{M} : B(\ell_{2}) \to B(\ell_{2})\| < +\infty.$$

If $M = [m_{ij}]$ is a Schur multiplier and the iterated row and
column limits
$$\lim_{j\to \infty} (\lim_{i\to \infty} m_{ij} ) = C \mbox{ and }
\lim_{i\to \infty} (\lim_{j\to \infty} m_{ij} ) = R$$
exist, then \cite{bennett} $C = R$. This shows in particular that
$$\left[ \frac{j-i}{i+j+1}\right]_{i,j\geq 1}$$
is not a Schur multiplier.

\begin{thm}\label{bennett}
	Let $(a_n)_{n\geq 1}$ be a sequence of reals which converges to $0$.
	Set
	$$b_n = \Delta^1(a_{n}) = a_n - a_{n+1}
	\mbox{ and } c_n = \Delta^2(a_{n}) = a_n - 2a_{n+1} + a_{n+2}.$$
	If the sequences $(a_n/n)$, $(b_n)$ and $(nc_n)$ are
	all absolutely summable, then the matrix
	$$\left[ \frac{(j-i)a_{i+j}}{i+j+1}\right]_{i,j\geq 1}$$
	is a Schur multiplier.
\end{thm}

\noi {\bf Proof.}
The proof will be based on the following criterion
due to Bennett \cite[Theorem 8.6]{bennett} : if $M = \left[
m_{ij}\right]_{i,j\geq 1}$
satisfies
$$\lim_{i}m_{ij} = \lim_{j}m_{ij} = 0$$
and
\begin{equation}
	\sum_{i,j=1}^{+\infty} \left| m_{i,j} - m_{i,j+1} - m_{i+1,j} +
	m_{i+1,j+1}\right| < + \infty ,
	\label{eq:bennett}
\end{equation}
then $M$ is a Schur multiplier.

 Let $m_{ij} = (j-i)a_{i+j}/(i+j+1)$. We have
$$\lim_{i}m_{ij} = \lim_{j}m_{ij} = 0.$$

In order to prove (\ref{eq:bennett}), we write
$$\sum_{i,j=1}^{\infty} \left| m_{i,j} - m_{i,j+1} - m_{i+1,j} +
	m_{i+1,j+1}\right| $$
$$= \sum_{n=1}^{\infty}\sum_{i+j=n} \left|
	\frac{(j-i)a_{i+j}}{i+j+1} - \frac{2(j-i)a_{i+j+1}}{i+j+2} +
	\frac{(j-i)a_{i+j+2}}{i+j+3}  \right|$$
$$\leq \sum_{n=1}^{\infty}n(n-1)\left| \frac{a_{n}}{n+1} -
\frac{2a_{n+1}}{n+2} + \frac{a_{n+2}}{n+3}  \right|$$
$$= \sum_{n=1}^{\infty}n(n-1)\left| \frac{c_{n}}{n+2} +
\frac{1}{n+2}\left[ \frac{a_{n}}{n+1} - 
\frac{a_{n+2}}{n+3}\right]\right|$$
$$= \sum_{n=1}^{\infty}n(n-1)\left| \frac{c_{n}}{n+2} +
\frac{b_{n}+b_{n+1}}{(n+1)(n+2)} + \frac{2a_{n+2}}{(n+1)(n+2)(n+3)} 
\right|$$
$$\leq \sum_{n\geq 1} n|c_{n}| + \sum_{n\geq 1} |b_{n}| + \sum_{n\geq 
1}
|b_{n+1}| + 2\sum_{n\geq 1}\frac{|a_{n+2}|}{n+2}.$$
All four sums are convergent by hypothesis. \myqed

\begin{cor}\label{cor:fe}
	Let $\ep > 0$. The matrices
$$E_{\ep} = \left[
\frac{j-i}{(i+j+1)(\log (i+j+1))^{1+\ep}}\right]_{i,j\geq 1}$$
	and
$$F_{\ep} = \left[ \frac{j-i}{(i+j+1)(\log (i+j+1))(\log
\log (i+j+1))^{1+\ep}}\right]_{i,j\geq 1}$$
	are Schur multipliers.
\end{cor}

\noi {\bf Proof.} Theorem \ref{bennett} applies for the sequences
$$a_{n} = \frac{1}{(\log n)^{1+\ep}}\, , n \geq 2$$
and
\begin{equation}
	a_{n} = \frac{1}{(\log n)[\log (\log n)]^{1+\ep}}\, , n \geq 2.
	\label{eq:32}
\end{equation}

We give the proof only for the second sequence (\ref{eq:32}).
Denote $r = 1+\ep$
and $\log_{2}(x) = \log (\log x)$. The
series
$$\sum_{n\geq 2}\frac{|a_{n}|}{n} = \sum_{n\geq 2}\frac{1}{n(\log
n)[\log_{2}(n)]^{r}}$$
is convergent since $r > 1$.

Consider the $C^1$ function
$$f(x) = \left(\log x\right)\left[ \log_{2}(x)\right]^{r}$$
with
$$f'(x) = \frac{(\log_{2}(x))^r + r(\log_{2}(x))^{r-1}}{x}.$$
We have
$$b_{n} = \frac{f(n+1) - f(n)}{\log (n)\left[
\log_{2}(n)\right]^r\log (n+1)\left[
\log_{2}(n+1)\right]^r}.$$
For each $n$ there is a point $\theta_{n}$ between $0$ and $1$ such
that
$$f(n+1) - f(n) = f'(n+\theta_{n}).$$
We obtain
$$\left| b_{n} \right| = \left| \frac{(\log_{2}(n+\theta_{n}))^r +
r(\log_{2}(n+\theta_{n}))^{r-1}}{(n+\theta_{n})\log (n)\left[
\log_{2}(n)\right]^r\log (n+1)\left[
\log_{2}(n+1)\right]^r}\right|$$
and thus
$$\left| b_{n} \right| \leq \frac{1+r}{n(\log n)[\log_{2}(n)]^r}$$
for sufficiently large $n$. Thus the sequence $(b_{n})$ is absolutely
summable.

Consider now the $C^2$ function $g(x) = 1/f(x)$, with its second
derivative given by
$$g''(x) = \frac{1}{x^2(\log x)^2(\log_{2}x)^{r}}$$
$$+ \frac{2}{x^2(\log
x)^3(\log_{2}x)^{r}} + \frac{3r}{x^2(\log x)^3(\log_{2}x)^{1+r}}$$
$$+ \frac{r}{x^2(\log
x)^2(\log_{2}x)^{1+r}}+ \frac{r(1+r)}{x^2(\log
x)^3(\log_{2}x)^{2+r}} .$$
For each n, there is $\eta_{n}$ between $0$ and $2$ such
that
$$c_{n} = g(n) - 2g(n+1) + g(n+2) = g''(n+\eta_{n}).$$
Using this representation of $c_{n}$ it can be proved that $(nc_{n})$
is absolutely summable.
\myqed

\section{A sufficient condition for similarity of $R(X)$
to $R(0)$}
The idea of the proof of the following theorem goes back to
\cite{foias/williams} and \cite{williams}.

\begin{thm}\label{thm:suff}
	Let $T_{2} \in B(H_{2})$ be an isometry ($T_{2}^*T_{2} = I_{H_{2}}$)
	and let $T_{1} \in B(H_{1})$ and $X : H_{1} \to H_{2}$ be bounded
	operators. Consider
	$$R(T_{2}^*,T_{1} ; X) = \left[\begin{array}{cc}
			T_{2}^* & X  \\
			0 & T_{1}
		\end{array}\right] .$$
	If
	\begin{equation}
		\sup_{n\geq 1}\| \sum_{j=0}^{n-1}T_{2}^{j+1}XT_{1}^j \| <
+\infty ,
		\label{eq:31}
	\end{equation}
	then $R(T_{2}^*,T_{1} ; X)$ is similar to
	$T_{2}^*\oplus T_{1} = R(T_{2}^*,T_{1} ; 0)$.
\end{thm}

\noi {\bf Proof.}  Let $\mathcal{L}$ be a
Banach limit \cite{conway}, that is a bounded
linear functional on $\ell_{\infty}(\C)$ such that $\1 =
\mathcal{L}(\1) = \|\mathcal{L}\|$ and $\mathcal{L}((x_{n+1})_{n\geq
0}) = \mathcal{L}((x_{n})_{n\geq 0})$ for every $(x_{n})_{n\geq 0} \in
\ell_{\infty}(\C)$. Here $\1 = (1,1, \ldots)$.

Consider the linear operator $Z : H_{1} \to H_{2}$ given by
$$\langle Zh_{1},h_{2}\rangle_{H_{2}} = \mathcal{L}\left( \langle
\sum_{j=0}^{n-1}T_{2}^{j+1}XT_{1}^j h_{1},h_{2}\rangle\right) .$$
Then (\ref{eq:31}) shows that $Z$ is well-defined and bounded.

We have
\begin{eqnarray*}
	\langle (T_{2}^*Z - ZT_{1})h_{1},h_{2}\rangle & = & \langle
	Zh_{1},T_{2}h_{2}\rangle - \langle ZT_{1}h_{1},h_{2}\rangle  \\
	 & = & \mathcal{L}\left( \langle \sum_{j=0}^{n-1}T_{2}^{j+1}XT_{1}^j
	 h_{1}, T_{2}h_{2}\rangle \right. \\
	 &  & -\left.\langle
	 \sum_{j=0}^{n-1}T_{2}^{j+1}XT_{1}^{j+1} h_{1},h_{2}\rangle\right)  
\\
	 & = & \mathcal{L}\left( \langle \sum_{j=0}^{n-1}T_{2}^{j}XT_{1}^j
	 h_{1}, h_{2}\rangle \right.\\
	 &  & - \left.\langle
	 \sum_{j=0}^{n-1}T_{2}^{j+1}XT_{1}^{j+1} h_{1},h_{2}\rangle\right)  
\\
	 & = & \langle Xh_{1},h_{2}\rangle - \mathcal{L}\left( \langle
	 T_{2}^nXT_{1}^nh_{1} , h_{2}\rangle\right).
\end{eqnarray*}
On the other hand,
\begin{eqnarray*}
	\mathcal{L}\left( \langle
	 T_{2}^nXT_{1}^nh_{1} , h_{2}\rangle\right) & = & \mathcal{L}\left(
\langle
	 T_{2}^nXT_{1}^nh_{1} , T_{2}^*T_{2}h_{2}\rangle\right)  \\
	 & = & \mathcal{L}\left( \langle
	 T_{2}^{n+1}XT_{1}^nh_{1} , T_{2}h_{2}\rangle\right)  \\
	 & = & \mathcal{L}\left( \langle
	 \sum_{j=0}^nT_{2}^{j+1}XT_{1}^jh_{1} , T_{2}h_{2}\rangle\right) \\
	 &  & -
	 \mathcal{L}\left( \langle
	 \sum_{j=0}^{n-1}T_{2}^{j+1}XT_{1}^jh_{1} ,
T_{2}h_{2}\rangle\right)  \\
	 & = & 0.
\end{eqnarray*}
Therefore $T_{2}^*Z - ZT_{1} = X$ and thus
$$
\left[\begin{array}{cc}
			I & -Z  \\
			0 & I
		\end{array}\right] \left[\begin{array}{cc}
			T_{2}^* & 0  \\
			0 & T_{1}
		\end{array}\right]\left[\begin{array}{cc}
			I & Z  \\
			0 & I
		\end{array}\right] = \left[\begin{array}{cc}
			T_{2}^* & X  \\
			0 & T_{1}
		\end{array}\right].
$$
We obtain that $R(T_{2}^*,T_{1} ; X)$ is similar to
$T_{2}^*\oplus T_{1} = R(T_{2}^*,T_{1} ; 0)$.
\myqed

\medskip

It follows from the above Theorem that if (\ref{eq:31}) holds and
$T_{1}$ is similar to a contraction, then
$R(T_{2}^*,T_{1} ; X)$ is similar to a contraction. 
We refer to
\cite{foias/williams,ccfw,clark,paulsen:yoneda,cassier} 
for related results. 

\begin{cor}
Let $T_{2} \in B(H_{2})$ be an isometry
	and let $T_{1} \in B(H_{1})$ and $X : H_{1} \to H_{2}$ be bounded
	operators such that the spectral radius of $T_{1}$ satisfies
	$r(T_{1}) < 1$. Then
	$$R(T_{2}^*,T_{1} ; X) = \left[\begin{array}{cc}
			T_{2}^* & X  \\
			0 & T_{1}
		\end{array}\right] $$
		is similar to a contraction.
\end{cor}

\noi {\bf Proof.} We have
$$\| \sum_{j=0}^{n-1}T_{2}^{j+1}XT_{1}^j \| \leq
\sum_{j=0}^{n-1}\|T_{2}^{j+1}\|\, \|X\|\, \|T_{1}^j\|$$
$$\leq \|X\| \sum_{j=0}^{+\infty}\|T_{1}^j\|$$
for each $n$. The last sum converges since $r(T_{1}) < 1$. By
Theorem \ref{thm:suff}, $R(T_{2}^*,T_{1} ; X)$ is similar to 
$T_{2}^*\oplus
T_{1}$. The operator $T_{2}^*$ is a contraction, while $T_{1}$ is 
similar
to a contraction by Rota's \cite{rota} theorem. Thus $T_{2}^*\oplus
T_{1}$ is similar to a contraction. \myqed

The following simple result characterizes when $R(T_{2}^*,T_{1} ; X)$
is power-bounded in terms of a condition related to (\ref{eq:31}) and
the power-bounded\-ness of $T_{1}$.
\begin{prop}
	Let $T_{2} \in B(H_{2})$ be an isometry ($T_{2}^*T_{2} = I_{H_{2}}$)
	and let $T_{1} \in B(H_{1})$ and $X : H_{1} \to H_{2}$ be bounded
	operators. Then
	$$R(T_{2}^*,T_{1} ; X) = \left[\begin{array}{cc}
			T_{2}^* & X  \\
			0 & T_{1}
		\end{array}\right] .$$
is power-bounded if and only if
	\begin{equation}
		\sup_{n\geq 1}
		\| T_{2}^{*n}\sum_{j=0}^{n-1}T_{2}^{j+1}XT_{1}^j \| < +\infty .
	\end{equation}
	and $T_{1}$ is power-bounded.
\end{prop}

\noi {\bf Proof.}  Set $R = R(T_{2}^*,T_{1} ; X)$. We have
$$R^n = \left[\begin{array}{cc}
			T_{2}^{*n} &
\sum_{j=0}^{n-1}T_{2}^{*(n-j-1)}XT_{1}^j  \\
			0 & T_{1}^n
		\end{array}\right] .
$$
Since $T_{2}$ is an isometry, we have
$$\sum_{j=0}^{n-1}T_{2}^{*(n-j-1)}XT_{1}^j = T_{2}^{*n}
\sum_{j=0}^{n-1}T_{2}^{j+1}XT_{1}^j$$
which gives the desired equivalence.
\myqed

\section{Proofs of the main results}
\noi {\bf Proof of Theorem \ref{thm:C}.}  By
Theorem \ref{thm:suff}, $R(S^{*(\infty)}, S^{(\infty)} ; X)$ is
similar to a contraction whenever
$$\sup_{n\geq 1}\|
\sum_{k=0}^{n-1}S^{(\infty)(k+1)}XS^{(\infty)k} \| < +\infty  .$$
The matrix of the operator
$S^{(\infty)k}$ is given by
$$\left[\begin{array}{cccc}
			0 & 0 & 0 & \cd  \\
			\vdots & \vdots & \vdots & \cd  \\
			0 & 0 & 0 & \cd  \\
			1 & 0 & 0 & \cd  \\
			0 & 1 & 0 & \cd  \\
			0 & 0 & 1 & \cd  \\
			\vdots & \vdots & \vdots & \ddots
		\end{array}\right] , $$
where the first non-zero entry occurs at position $(k,0)$. The matrix
of the operator $S^{(\infty)(k+1)}XS^{(\infty)k}$ is
$$\left[\begin{array}{cccc}
			0 & 0 & 0 & \cd  \\
			\vdots & \vdots & \vdots & \cd  \\
			0 & 0 & 0 & \cd  \\
			X_{0k} & X_{0,k+1} & X_{0,k+2} & \cd  \\
			X_{1,k} & X_{1,k+1} & X_{1,k+2} & \cd  \\
			\vdots & \vdots & \vdots & \ddots
		\end{array}\right] , $$
where the first non-zero entry occurs at position $(k+1,0)$. This
gives that the entries of the matrix $\tilde{A}^{(n)}(X)$ of the 
operator
$$\sum_{k=0}^{n-1}S^{(\infty)(k+1)}XS^{(\infty)k}$$
are $0$ on the first
line and
$$\tilde{A}^{(n)}(X)_{ij} = X_{i-1,j} + X_{i-2,j+1} + \ldots
	X_{i-1-\min(i-1,n-1),j+\min(i-1,n-1)} $$
for $i\geq 1$. Then the operator represented by
$\tilde{A}^{n}(X)$ is bounded whenever
the operator represented by $A^{(n)}(X)$ is.
\myqed

\medskip

\noi {\bf Proof of Theorem \ref{thm:charact}.} By 
\cite{foias/williams}, 
the Foguel-Hankel operator
$$R(\Gamma) = R(S^{*(\infty)}, S^{(\infty)} ; \Gamma)$$
is similar to a
contraction if and only if there is a bounded solution $Y$ of the 
equation
\begin{equation}
	S^{*(\infty)}Y - YS^{(\infty)} = \Gamma,
	\label{eq:sylv}
\end{equation}
if and only if $R(\Gamma)$ is similar to $S^{*(\infty)}\oplus
S^{(\infty)}.$
One completes the proof by observing that each of the 
three operators, $-\Gamma D, D^*\Gamma$ and $(-\Gamma D 
+D^*\Gamma)/2$ 
is a formal solution to this commutator equation.
\myqed

\medskip

In the case of a scalar symbol $\Gamma_{f},$ one can say considerably 
more.
If a matrix $Y$ is any formal solution(bounded or unbounded) of 
(\ref{eq:sylv}), then $-Y^{t}$
is another solution of the same(possibly infinite)
norm. Here $Y^{t}$ is the transpose of $Y$. Indeed, $(\Gamma_f)^t =
\Gamma_f$ and $S^{*} = S^{t}$.
Hence $[Y - Y^{t}]/2$ is a solution of no greater norm.
Now if $Y$ is a solution of (\ref{eq:sylv}),
then any other solution differs from $Y$ by a matrix $G$ satisfying
$$S^{*}G = GS,$$
that is by a formal Hankel
matrix $G$. But
$$(Y+G)- (Y+G)^{t} = Y-Y^{t}$$
and so
if $Y$ is any solution of (\ref{eq:sylv}), then
$[Y-Y^{t}]/2$ is the solution of minimum norm.

Now
$$Y_{0} = -\Gamma_f D = \left[ -ja_{i+j-1}\right]_{i,j\geq 0}$$
is a formal solution of (\ref{eq:sylv}) and thus
$$([ \frac{(j-i)a_{i+j-1}}{2}]_{i,j\geq 0})^t 
= \frac{(\Gamma_{f} D - D^*\Gamma_{f})^{t}}{2} =
 \frac{[Y_{0}-Y_{0}^{t}]}{2}$$
is the solution of minimum norm.

Thus, we see immediately as an application of \cite{foias/williams} 
that $R(\Gamma_f)$ is similar to a contraction if and only 
if $\Gamma_f D -D^*\Gamma_f$ is bounded.

Unfortunately, the above arguement fails in the operator-valued 
Hankel setting.  The difficulty is that the transpose of a bounded 
operator matrix need not be bounded. This is essentially because the 
transpose map is not completely bounded.  For this reason it is not 
known what the minimum norm solution of (\ref{eq:sylv})  is in the 
operator case. This makes it difficult to obtain necessary and 
sufficient conditions for the existence of bounded solutions.
 
\medskip

\noi {\bf Proof of Theorem \ref{thm:A}.} 
By Theorem \ref{thm:charact} it is sufficient to show that
$$[(j-i)\alpha_{i+j-1}C_{i+j-1}]_{i,j\geq 1}$$
is bounded. Indeed, this matrix is obtained from the matrix 
$$[(j-i)\alpha_{i+j-1}C_{i+j-1}]_{i,j\geq 0}$$
of Theorem \ref{thm:charact} by adding 
a bounded row and column.

Let $\ep > 0$. The matrix
$[(j-i)\alpha_{i+j-1}C_{i+j-1}]_{i,j\geq 1}$ can be viewed as the
Schur product of $F_{\ep/2}$ and
$$\left[(i+j+1)\left(\log (i+j+1)\right)\left(\log (\log (i+j+1))
\right)^{(\ep+2)/2}\alpha_{i+j-1}C_{i+j-1}\right]_{i,j\geq 1} .$$

The last matrix represents a bounded Hankel operator.
Indeed, we have \cite{davidson/paulsen}
$$\left|\!\left|\left[(i+j+1)\left(\log (i+j+1)\right)\left(\log (\log
(i+j+1))
\right)^{(\ep+2)/2}\alpha_{i+j-1}C_{i+j-1}\right]\right|\!\right|$$
$$= \left|\!\left|\left[(i+j+1)\left(\log (i+j+1)\right)
\left(\log (\log (i+j+1))
\right)^{(\ep+2)/2}\alpha_{i+j-1}e_{i+j-1}\right]\right|\!\right|$$
$$= \sum_{k\geq 1} (k+2)^2\left(\log (k+2)\right)^2\left[ \log (\log
	(k+2)\right]^{2+\ep}|\alpha_{k}|^2$$
which is convergent.

Since $F_{\ep/2}$ is a Schur multiplier by Corollary \ref{cor:fe}, 
their Schur product gives 
\cite[Theorem 4.1]{davidson/paulsen} a bounded operator. \myqed

\end{document}